\documentclass[a4paper,11pt]{article}
\usepackage[utf8x]{inputenc}

\usepackage[title]{appendix}

\usepackage{a4wide}
\usepackage{amssymb}
\usepackage{amsmath}
\usepackage{amsthm}
\usepackage{mathtools}
\usepackage[mathscr]{euscript}
\usepackage[normalem]{ulem}
\usepackage[colorlinks=true]{hyperref}
\usepackage[usenames,dvipsnames]{xcolor}
\usepackage{appendix}
\usepackage{physics}

\usepackage{graphicx}

\numberwithin{equation}{section}

\theoremstyle{plain}

\newtheorem{Definition}{Definition}[section]

\newtheorem{Theorem}[Definition]{Theorem}
\newtheorem{Proposition}[Definition]{Proposition}
\newtheorem{Lemma}[Definition]{Lemma}
\newtheorem{Corollary}[Definition]{Corollary}
\newtheorem{Conjecture}[Definition]{Conjecture}

\theoremstyle{definition}

\newcommand{\R}{\mathbb R}

\usepackage[xcolor]{changebar}
\cbcolor{blue}
\newcommand{\Ric}{\mathrm{Ric}}

\usepackage[inline]{enumitem}
\newcommand{\enumlabelformat}{\roman}
\newcommand{\enumlabelfont}[1]{#1}
\newlength{\thelabelsep}
\setlist{labelsep=\thelabelsep}
\setlist[enumerate,1]{font=\enumlabelfont,label=(\enumlabelformat*),leftmargin=2.5em}
\setlist[itemize]{leftmargin=2.5em,label=$-$}

\newcounter{inlineenum}
\renewcommand{\theinlineenum}{\enumlabelformat{inlineenum}}

\let\epsilon\varepsilon
\let\phi\varphi

\makeatletter
\let\save@mathaccent\mathaccent
\newcommand*\if@single[3]{%
  \setbox0\hbox{${\mathaccent"0362{#1}}^H$}%
  \setbox2\hbox{${\mathaccent"0362{\kern0pt#1}}^H$}%
  \ifdim\ht0=\ht2 #3\else #2\fi
  }
\newcommand*\rel@kern[1]{\kern#1\dimexpr\macc@kerna}
\newcommand*\widebar[1]{\@ifnextchar^{{\wide@bar{#1}{0}}}{\wide@bar{#1}{1}}}
\newcommand*\wide@bar[2]{\if@single{#1}{\wide@bar@{#1}{#2}{1}}{\wide@bar@{#1}{#2}{2}}}
\newcommand*\wide@bar@[3]{%
  \begingroup
  \def\mathaccent##1##2{%
    \let\mathaccent\save@mathaccent
    \if#32 \let\macc@nucleus\first@char \fi
    \setbox\z@\hbox{$\macc@style{\macc@nucleus}_{}$}%
    \setbox\tw@\hbox{$\macc@style{\macc@nucleus}{}_{}$}%
    \dimen@\wd\tw@
    \advance\dimen@-\wd\z@
    \divide\dimen@ 3
    \@tempdima\wd\tw@
    \advance\@tempdima-\scriptspace
    \divide\@tempdima 10
    \advance\dimen@-\@tempdima
    \ifdim\dimen@>\z@ \dimen@0pt\fi
    \rel@kern{0.6}\kern-\dimen@
    \if#31
      \overline{\rel@kern{-0.6}\kern\dimen@\macc@nucleus\rel@kern{0.4}\kern\dimen@}%
      \advance\dimen@0.4\dimexpr\macc@kerna
      \let\final@kern#2%
      \ifdim\dimen@<\z@ \let\final@kern1\fi
      \if\final@kern1 \kern-\dimen@\fi
    \else
      \overline{\rel@kern{-0.6}\kern\dimen@#1}%
    \fi
  }%
  \macc@depth\@ne
  \let\math@bgroup\@empty \let\math@egroup\macc@set@skewchar
  \mathsurround\z@ \frozen@everymath{\mathgroup\macc@group\relax}%
  \macc@set@skewchar\relax
  \let\mathaccentV\macc@nested@a
  \if#31
    \macc@nested@a\relax111{#1}%
  \else
    \def\gobble@till@marker##1\endmarker{}%
    \futurelet\first@char\gobble@till@marker#1\endmarker
    \ifcat\noexpand\first@char A\else
      \def\first@char{}%
    \fi
    \macc@nested@a\relax111{\first@char}%
  \fi
  \endgroup
}
\makeatother

\newcommand\beq{\begin{equation}}
\newcommand\eeq{\end{equation}}
\newcommand\ben{\begin{enumerate}}
\newcommand\een{\end{enumerate}}
\newcommand\bit{\begin{itemize}}
\newcommand\eit{\end{itemize}}
\newcommand{\p}{\partial}
\renewcommand{\S}{\Sigma}

\renewcommand{\d}{\nabla}

\newcommand\e{\varepsilon}

\newcounter{mnotecount}

\let\div\relax
\DeclareMathOperator{\div}{div}

\title{Cosmological rigidity in $3$ spacetime dimensions}

\author{Gregory J.\ Galloway\thanks{Department of Mathematics, University of Miami, Coral Gables, Florida, USA (\textsf{galloway@math.miami.edu})}\\Eric Ling\thanks{Department of Mathematics, Florida State University, Tallahassee, Florida, USA (\textsf{eling@fsu.edu})}\\Argam Ohanyan\thanks{Department of Mathematics, University of Toronto, Ontario, Canada (\textsf{argam.ohanyan@utoronto.ca})}}

\begin{document}

\date{\today}


\maketitle

\begin{abstract} We establish rigidity theorems for cosmological spacetimes in $3$ spacetime dimensions. In the setting of the Bartnik Splitting Conjecture, under the additional assumptions of null geodesic completeness and the dominant energy condition, we show that the spatial topology is either spherical or the spacetime splits as a product. Our methods rely on the theory of marginally outer trapped surfaces, known rigidity results in the presence of a timelike or null line,  and more recent rigidity results in the context of Penrose's singularity theorem due to the first two authors, from which the existence of a foliation by totally geodesic null hypersurfaces can be concluded. Our second main result closes the argument by establishing a splitting theorem for cosmological and null geodesically complete $3$-dimensional spacetimes which satisfy the dominant energy condition and which admit a foliation of the aforementioned type. In its proof, we make use of the local Kundt structure of the metric and construct a global parallel generator for the null foliation.

\vspace{1em}

\noindent
\emph{Keywords:} Bartnik's splitting conjecture, splitting, cosmological spacetime, null foliation
\medskip

\noindent
\emph{MSC2020:} 53C50, 53C24

\end{abstract}

\tableofcontents

\section{Introduction}\label{section: introduction}

In the early 1980's,  S.T.\ Yau put forth the problem of establishing the {\it rigidity} of the Hawking-Penrose singularity theorems.  From the standpoint of rigidity one would like to eliminate from the singularity theorems conditions like the so-called {\it generic condition} and retain only weak curvature inequalities such as the strong energy condition (SEC), $\Ric(X,X) \ge 0$ for all timelike vectors $X$.  This viewpoint became formalized in a 1988 conjecture of R. Bartnik \cite{Bart88}   concerning spatially closed spacetimes.  

\begin{Conjecture}[Bartnik splitting conjecture]
\label{conj}
Let $(M,g)$ be a  globally hyperbolic spacetime with compact Cauchy surfaces which obeys the SEC. If $(M,g)$ is timelike geodesically complete then $(M,g)$ splits, i.e.\ $(M, g)$ is isometric to the product
$(\R\times V,-dt^2 + h)$, where $(V,h)$ is a compact Riemannian manifold with nonnegative Ricci curvature.
 \end{Conjecture}

Thus, according to the conjecture, a spatially closed spacetime obeying the SEC can avoid being  `singular' (i.e. timelike geodesically incomplete) only under very special circumstances.
Despite various attempts, and partial results, the conjecture remains open. We refer to the survey \cite{Gallowaysurvey} as well as the newer works \cite{GallowayLingExpandingCMC, GMO, ohanyan2026cosmological, MinguzziWall} for more details on the aforementioned developments.  

As an approach to such rigidity questions, Yau, in his famous 1982 problem section \cite{Yau}, posed the problem of establishing a Lorentzian analogue of the Cheeger-Gromoll splitting theorem of Riemannian geometry.  This was settled in a series of papers by the end of the 1980's \cite{Eschenburg1988, GalLST, Newman} (see also  \cite{BGMOS2026} for a recent proof of elliptic flavor). The most basic version, due to Eschenburg, is as follows.  

\begin{Theorem}[The Lorentzian Splitting theorem] 
\label{Theorem: Lorentzian splitting}
{\it Let $(M,g)$ be a globally hyperbolic, timelike geodesically complete spacetime, satisfying the SEC.
If $M$ admits a timelike line, then $M$ {\it splits}, i.e.\
$(M, g)$ is isometric to the product $(\R\times V,-dt^2 + h)$  where $(V,h)$ is a  geodesically complete Cauchy surface with nonnegative Ricci curvature. }
\end{Theorem}

A timelike line is an inextendible timelike geodesic each segment of which is maximal; a null line is defined analogously (see Section \ref{section: prelim}). The Bartnik splitting conjecture then follows if one can establish the existence of a timelike line.  As discussed in Section \ref{section: 3d Bartnik}, a globally hyperbolic spacetime with compact Cauchy surfaces always admits  a timelike line or a null line.  While there are known conditions under the assumptions of the Bartnik Splitting Conjecture (BSC) that guarantee the existence of a timelike line, the possibility of a null line has not been ruled out under those assumptions alone. In \cite{EhrlichGalloway} a spacetime was constructed which is causal (timelike and null) geodesically complete, with compact Cauchy surfaces, that does not contain any timelike lines.  However this spacetime does not satisfy the SEC.  

The possibility of a null line under the assumptions of the BSC motivated the first-named author to study the rigidity of spacetimes admitting a null line; see \cite{GalNullMaxPrin} and Theorem \ref{thm:nullsplit} in the next section.  This so-called `null splitting theorem' requires null geodesic completeness.  This suggests  considering the BSC under the assumption of causal geodesic completeness, rather than just timelike geodesic completeness.  The BSC remains open even under this stronger completeness assumption.  

In Section \ref{section: 3d Bartnik} we present a Bartnik-type splitting result for three-dimensional spacetimes that are assumed to be causal geodesically complete.  Various elements arise in its proof. One of the key ingredients in the proof is a result from \cite{GallowayLing2025}, which implies
that for the three-dimensional spacetimes being considered, under certain circumstances, they admit a foliation by totally geodesic null hypersurfaces.  In Section \ref{sec: null foliation} we present a splitting result for a class of three-dimensional spacetimes that admit such a foliation.  This, in turn,  is used to complete the proof of our Bartnik-type splitting result.
In the next section we recall some basic notions from spacetime geometry, and review some properties of null hypersurfaces used in Sections \ref{section: 3d Bartnik} and \ref{sec: null foliation}, and of marginally outer trapped surfaces used  in Section \ref{section: 3d Bartnik}.

\section{Preliminaries}
\label{section: prelim}

In this section we review some basic notions and results in spacetime geometry needed in the subsequent sections.  By a spacetime we mean a smooth, connected, time-oriented Lorentzian manifold 
$(M,g)$.  We adopt the signature convention $(-+\cdots+)$.  
We recall a couple of basic causal theoretic notions. 
A spacetime $(M,g)$ is globally hyperbolic provided (i) $M$ is  causal and  (ii)  the `causal diamonds' $J^+(p) \cap J^-(q)$, $p,q \in M$ are compact.  
A Cauchy surface in $M$ is a subset $V$ of $M$ that is met exactly once by every inextendible causal curve in $M$. Cauchy surfaces are in general $C^0$ acausal hypersurfaces, and any two Cauchy surfaces in $M$ are homeomorphic.   

It is a fundamental fact that $(M,g)$ is globally hyperbolic if and only if it admits a smooth spacelike Cauchy surface.  In fact,  $(M,g)$ is globally hyperbolic if and only if it admits a  Cauchy temporal function, i.e. a smooth time function with past timelike gradient,
$\tau:M\to\mathbb R$,
such that its level sets $V_t:=\tau^{-1}(t)$
are smooth spacelike Cauchy hypersurfaces (\cite{BernalSanchez2005}). Henceforth, for the purposes of this paper, we will only consider Cauchy surfaces that are smooth spacelike hypersurfaces. 

In our results we will be using several energy conditions:

\vspace{-.05in}
\ben
\item[(i)] Strong energy condition (SEC):  $\Ric(X,X) \ge 0$ for all timelike vectors $X$.

\vspace{-.05in}
\item[(ii)] Null energy condition (NEC): $\Ric(X,X) \ge 0$ for all null vectors $X$

\vspace{-.05in}
\item[(iii)] Dominant energy condition (DEC):  $G(X,Y) \ge 0$ for all future directed causal vectors $X,Y$, where $G$ is the Einstein tensor,
\begin{equation}
    G=\Ric - \frac{1}{2} Rg.
\end{equation}
\een
Here we are using the following convention for the curvature tensor,
\begin{equation}
    R(X,Y)Z = \nabla_X \nabla_Y Z - \nabla_Y \nabla_X Z - \nabla_{[X,Y]}Z ,
\end{equation}
where $\d$ is the Levi-Civita connection.  Note that the SEC and DEC imply the NEC, and when the scalar curvature vanishes the DEC implies the SEC.

\subsection{Null hypersurfaces and null foliations}
\label{subsec: prelim null}

We now recall some aspects of null hypersurfaces.  A codimension one linear subspace $H \subseteq T_xM$ is called \emph{null} if $g_x|_{H^2}$ is degenerate. It is easily seen that $H$ is degenerate if and only if $H^\perp \subseteq H$, and $\dim H^\perp = 1$.

\begin{Definition}
A smooth immersed hypersurface $N$ in $(M,g)$ is called
null if $T_pN$ is a null hyperplane for every $p\in N$.
Its `radical' is the smooth line bundle over $N$
\[
\operatorname{Rad}(TN)
:=TN\cap(TN)^\perp=(TN)^\perp.
\]
A smooth nowhere-zero section $k$ of this line bundle
is called a null generator field.
\end{Definition}
 Unless otherwise stated, we shall take null generator fields to be future pointing.  It is always possible to choose such a  null generator field $k$. Further, all other null generators are of the form $f k$ with $f > 0$ smooth. At each point of $N$, $k$ is the direction of degeneracy
 and $k^\perp = TN$.  A null generator field $k$ for $N$ may be viewed as a null vector field on $N$.  

The null second fundamental form associated to a null generator field $k$ is the covariant $2$-tensor $\Psi_k$ given by, $\Psi_k(X,Y) = g(\d_Xk,Y)$, for $X, Y \in\Gamma(TN)$.  By standard arguments, $\Psi_k$ is symmetric in $X$ and $Y$. 
 It is easily seen that if $\widetilde{k} = f k$, then $\Psi_{\widetilde{k}} = f \Psi_k$.  A smooth null hypersurface $N$ is said to be {\it totally geodesic} if  $\Psi_k = 0$ for some null generator field $k$.  (Note that this is well-defined.) 

We have the following  equivalent version of a null hypersurface being totally geodesic. 

\begin{Lemma}
\label{Lemma: nabla X k in span k}
Let $N$ be a smooth null hypersurface and let $k$ be a null generator field.
 Then $N$ is totally geodesic if and only if
\[
\nabla_XY\in\Gamma(TN)
\]
for all smooth vector fields $X,Y$ tangent to $N$.
\end{Lemma}

\begin{proof} For vector fields $X,Y$ tangent to $N$, 
the equation $X(g(k,Y))=0$ gives
\[
g(\nabla_Xk,Y)=-g(k,\nabla_XY).
\]
If $N$ is totally geodesic, the left-hand side vanishes.
Thus $\nabla_XY  \in k^\perp = TN$.

Conversely, if for all $X, Y \in \Gamma(TN)$, $\nabla_XY \in  \Gamma(TN)$, then the right-hand side vanishes and $N$ is totally geodesic.
\end{proof}

\smallskip
As a consequence of this lemma, for $N$ a totally geodesic null hypersurface and $X \in \Gamma(TN)$, 
$$
\d_Xk = \lambda k
$$  
for some smooth function $\lambda$ on $N$ (in general depending on $X$).  Setting $X = k$, we then have $\d_kk = \lambda k$ for some  $\lambda$, i.e. the integral curves of $k$ are (in general) pregeodesics.  (We should point out that this fact holds even if $N$ is not totally geodesic.)  

\smallskip
The following sections will involve spacetimes that admit {\it foliations} by totally geodesic null hypersurfaces.  

\begin{Definition}
A smooth codimension-one foliation $\mathcal F$ of $M$
is a partition into connected immersed hypersurfaces,
called leaves, such that locally there are smooth
coordinates $(u,v,x)$ in which the connected pieces
of the leaves are given by $u=\mathrm{constant}$.

The foliation is called null if every leaf is a null
hypersurface, and totally geodesic if every leaf is
totally geodesic.
\end{Definition}

For a smooth null foliation $\mathcal F$, its radical line bundle $L:=\mathcal{F}^\perp$ is a line bundle over $M$. It is always possible to make a choice for future null generator $k$ for $L$ (i.e., a nonzero section of $L$).  We may view such a null generator as a null vector field on $M$ tangent to the leaves of  $\mathcal F$.

\begin{Lemma}[Complete affine generator]
\label{Lemma: complete affine generator}
Let $(M,g)$ be a globally hyperbolic,
null geodesically complete spacetime and let $\mathcal F$
be a smooth foliation by totally geodesic null hypersurfaces.
Then its radical line bundle admits a smooth future-directed null
generator $\widehat k$ such that
\[
\nabla_{\widehat k}\widehat k=0.
\]
Moreover, $\widehat k$ is complete.
\end{Lemma}

\begin{proof}
By applying Lemma \ref{Lemma: nabla X k in span k} as above to the leaves of the foliation we obtain a smooth future-directed null generator  $k$ for  $\mathcal F$ satisfying $\nabla_k k=\kappa k$, for some smooth function $\kappa$ on $M$.
Fix a smooth spacelike Cauchy hypersurface $V$.
Each maximal integral curve of $k$, regarded as an
unparametrized causal curve, is inextendible and hence
meets $V$ exactly once.

Solve
\[
k(f)+\kappa f=0,\qquad f|_V=1,
\]
along these curves. The solution is positive.
Since $k$ is transverse to $V$, smooth dependence
of its flow and of the scalar ODE gives $f\in C^\infty(M)$.

Setting $\widehat k=fk$, we obtain
\[
\nabla_{\widehat k}\widehat k
=f\bigl(k(f)+\kappa f\bigr)k=0.
\]
Its integral curves are therefore affinely parametrized
null geodesics. Null geodesic completeness and uniqueness 
of the geodesic equation imply that these integral curves
extend to all of~$\mathbb R$.
\end{proof}

Spacetimes foliated by totally geodesic null hypersurfaces arise in a number of situations, 
see e.g.\ \cite{GalGrafLing, GallowayLing2025}; aspects of  \cite{GallowayLing2025} will be considered in the next section. Both papers begin by establishing the existence of a single totally geodesic null hypersurface.  This is done by applying the  `null splitting theorem' \cite{GalNullMaxPrin}, which requires the existence of a {\it null line}.  

A \emph{causal line} in spacetime is an inextendible causal geodesic each segment of  which maximizes Lorentzian distance.  A causal line can be either timelike or null. Null lines can also be characterized as inextendible causal geodesics that are globally achronal. 

\begin{Theorem}[{\cite[Thm. IV.1]{GalNullMaxPrin}}]
\label{thm:nullsplit}
Let $(M,g)$ be a null geodesically complete spacetime which obeys the NEC and contains a null line $\eta$. Then $\eta$ is contained in a smooth closed achronal totally geodesic null hypersurface.
\end{Theorem} 

\subsection{Marginally outer trapped surfaces}
\label{subsec: prelim MOTS}

As observed in the next section, there is a connection between totally geodesic null hypersurfaces and marginally outer trapped surfaces. Here we discuss some basic properties of the latter that will be needed there. 

Let $(M,g)$ be a spacetime of dimension $n\ge 3$,  and let $V$  be a smooth spacelike hypersurface in $M$ with induced metric $h$ and second fundamental form $K$. Let 
$\Sigma$ be a connected closed (compact without boundary) two-sided hypersurface in $V$.  Then $\Sigma$ admits a smooth unit normal field $\nu$ tangent to $V$; by convention we refer to 
$\nu$ as outward pointing.  Let $u$ be the future directed unit normal field to $V$. $\S$ admits two future directed null normal vector fields, $\ell^{\pm} = u \pm \nu$.  The \emph{null second fundamental forms} $\chi^{\pm}$ are given by 
\[
\chi^{\pm}(X,Y) \,=\,  g(\nabla_X \ell^{\pm}, Y) 
\]
for all $X,Y \in T_p\S$ and $p \in \S$.  Tracing the null second fundamental forms with respect to the induced metric on $\Sigma$, we obtain the \emph{null expansion scalars} 
$\theta^{\pm}$, 
\[
\theta^{\pm} \,=\, \rm{tr}_\Sigma \chi^{\pm} \,=\, \div_\S \ell^{\pm} \,=\, \rm{\tr}_\S K \pm H,
\]
where $H = \div_\S \nu$ is the mean curvature of $\Sigma$ within $V$.

$\S$ is \emph{outer trapped} if $\theta^+ < 0$, \emph{weakly outer trapped} if $\theta^+ \leq 0$, and \emph{marginally outer trapped} if $\theta^+ = 0$. In the latter case, $\S$ is referred to as a \emph{marginally outer trapped surface} (or MOTS for short).  (We remark that the \emph{sign} of 
$\theta^+$ does not depend on the particular scaling of the null vector fields $\ell^{\pm}$.)
Note that  in the case  $V$ is totally geodesic ($K = 0$), a MOTS is  a minimal surface in $V$. 

Let $\S \subset V$ be a MOTS, and let $\{\Sigma_t\}_{|t| < \e}$ be a variation of $\S$ given by 
$\S_t \colon x \mapsto \exp_x\big(t \phi(x)\nu(x)\big)$ for some $\phi \in C^\infty(\Sigma)$. Let $\theta^+(t)$ denote the outward null expansion of $\S_t$. From e.g. \cite[Eq.~1]{AMSPRL} (see also \cite[Prop.\ 7.32]{DLee}), we have that
\[
\frac{\p \theta^+}{\p t}\bigg|_{t=0} \,=\, L(\phi), 
\]
where
\[
L(\phi)\,=\,-\Delta \phi + 2 g(X, \nabla \phi) + (Q + \div X - |X|^2)\phi.
\]
Here $X$ is the vector field on $\Sigma$ metrically dual to the one-form $K(\nu, \cdot)|_\Sigma$ and
\[
Q \,=\, \frac{1}{2}S - G(u, \ell^+) - \frac{1}{2}|\chi^+|^2,
\]
where $S$ is the scalar curvature on $\S$. 

When $V$ is totally geodesic ($K = 0$),
$\theta^+ = H$, $X = 0$, and $L$ reduces to the classical stability operator of minimal surface theory. In analogy with the minimal surface case, we refer to $L$ as the \emph{MOTS stability operator}. Note, however, that $L$ is not in general self-adjoint. 

Nevertheless, it is shown in \cite{AMSPRL,AMS08} (see also 
\cite[Thm.~A.10]{DLee}) that, as a consequence of the Krein-Rutman theorem, $L$ admits a {\it principal eigenvalue} $\lambda_1$, which is real and simple. Furthermore there exists an associated principal eigenfunction $\phi$ which can be chosen to be strictly positive.  
In analogy with the minimal surface case, we say a MOTS $\S$ is \emph{stable} if $\lambda_1 \geq 0$.

The following variant of the Penrose singularity theorem (see e.g. \cite[Thm. 7.1]{AMMS}) will be used to establish the stability of certain marginally outer trapped surfaces arising in the next section.

\begin{Theorem}[One-sided Penrose]\label{thm: one-sided Penrose}
Let $(M,g)$ be a globally hyperbolic spacetime satisfying the null energy condition with a spacelike Cauchy hypersurface $V$. Let $\Sigma$ be a connected, closed hypersurface in $V$ which separates $V$. If $\S$ is outer trapped (i.e., $\theta^+ < 0$) with respect to a component of $V \setminus \S$ whose closure in $V$ is noncompact, then $M$ is future null geodesically incomplete.
\end{Theorem}

\section{A $2+1$ Bartnik-type splitting}
\label{section: 3d Bartnik}

Our main Bartnik-type result is the following.

\begin{Theorem}\label{Theorem: main}
Let $(M,g)$ be a null geodesically complete, three-dimensional spacetime satisfying the dominant energy condition. Let $V$ be a smooth compact spacelike Cauchy surface in $M$. Then either $V$ is the two-sphere/projective plane, or $(M,g)$ contains no null lines.

Hence, in the latter case, if in addition $(M,g)$ satisfies the strong energy condition and is timelike geodesically complete, then $(M,g)$ splits isometrically as $(\R \times V, \, -dt^2 + h)$. Moreover, if $V$ is neither the two-sphere nor the projective plane, then the metric $h$ is flat.
\end{Theorem}

\noindent \emph{Remark.} We note that de Sitter spacetime is null geodesically complete, satisfies the dominant energy condition, has compact Cauchy surfaces which are topologically spheres, and contains null lines. Thus, under the initial assumptions of Theorem \ref{Theorem: main}, the spherical case cannot be removed.

\begin{proof}[Proof of Theorem \ref{Theorem: main}]
We suppose that $V$ is neither the two-sphere nor the projective plane and that $(M,g)$ contains a null line $\gamma$. We will lead these assumptions to a contradiction. Since null lines lift to covers, we can assume that $V$ is oriented by passing to the orientable double cover. Theorem \ref{thm:nullsplit} implies that $\gamma$ is contained in a smooth, closed, achronal, totally geodesic null hypersurface $N$, which we take to be connected. The intersection $\Sigma=N\cap V$ is a compact connected embedded one-dimensional submanifold of $V$. Thus $\Sigma$ is  diffeomorphic to a circle. Projecting a global future-directed null generator field $k$ of $N$ onto $TV$ along $\Sigma$, and normalizing, gives a global unit normal $\nu$ to $\Sigma$ in $V$. Consequently, $\Sigma$ is two-sided. Taking $\nu$ as the outward normal, the total geodesicity of $N$ ($\Psi_k = 0$) implies that $\Sigma$ is a MOTS, i.e., $\theta^+ = 0$  (in fact, $\chi^+ = 0$).

Assume first that $\Sigma$ separates $V$.  Since $V$ is not the two-sphere, at least one component of $V \setminus \Sigma$, has closure which is not a disc. Reversing
the time orientation if necessary, we may assume that $\nu$
points into said component. Let $\lambda_1$ be the principal eigenvalue for the MOTS stability operator with respect to $\Sigma$. We will show that $\lambda_1 = 0$.  Let $\phi$ be a corresponding positive eigenfunction. Then
\[
    \frac{d\theta^+}{dt} = L(\phi) = \lambda_1\phi = - \Delta \phi + 2\langle X, \nabla \phi\rangle + (-P + \div X - |X|^2)\phi,
\]
where $P \geq 0$ by the dominant energy condition. Dividing by $\phi$ and setting $Y = X - \nabla \log \phi$, we find $\lambda_1 = \div Y - |Y|^2 - P$ (cf.\ \cite{GalSchoen2006}).  Thus $\lambda_1|\Sigma| = - \int_\Sigma (|Y|^2 + P)$, where $|\Sigma|$ denote the length of $\Sigma$. Therefore $\lambda_1 \leq 0$.

Suppose $\lambda_1<0$. Varying $\Sigma$ outward with initial
velocity $\phi\nu$, 
produces a nearby embedded circle $\Sigma_1$, with outward unit normal $\nu_1$, isotopic to $\Sigma$
with $\theta_1^+<0$. Let $A$ and $B$ denote the closures of the two components of $V \setminus \Sigma_1$, so that $V = A \cup_{\Sigma_1} B$. Assume $B$ is the closure of the side into which $\nu_1$ points. Cap $B$ with a disc, $\widehat{B} = D^2 \cup_{\Sigma_1} B$. Then $\widehat{B}$ is a closed orientable surface of positive genus, so its universal cover is $\R^2$. Remove all the lifted capping discs. The resulting surface, $B'$ has infinitely many boundary circles, and the restriction $B' \to B$ is a covering such that each boundary circle maps diffeomorphically onto $\Sigma_1$. Attach a copy of $A$ to each boundary circle. The resulting surface $\widetilde{V}$ is the desired covering space. (The subgroup of $\pi_1(V)$ corresponding to $\widetilde{V}$ is the normal closure of $\iota_* \pi_1(A)$, where $\iota \colon A \hookrightarrow V$ is the inclusion map.) The covering
$p\colon\widetilde V\to V$  induces a spacetime covering
$P\colon\widetilde M\to M$ for which $\widetilde V$ is a Cauchy
surface for $\widetilde{M}$. The maps $p$ and $P$ are related via a commutative diagram, as described at the beginning of Section 3.2 in
\cite{GallowayLing2022}. Each component of the lift of $\Sigma_1$ maps diffeomorphically onto $\Sigma_1$ and separates $\widetilde{V}$ into two components, the outer one having noncompact closure. Let $\widetilde{\Sigma}_1$ denote one such component of the lift. Its outer expansion is negative, $\widetilde{\theta}_1^+ < 0$. 
Theorem \ref{thm: one-sided Penrose} therefore implies that $\widetilde{M}$ is future null geodesically incomplete. This contradicts null geodesic completeness of $M$, which lifts to $\widetilde{M}$. Thus $\lambda_1 = 0$.

Since $\lambda_1 = 0$, by Lemma 1 in \cite{GallowayLing2025}, there is a neighborhood $U \approx (-\varepsilon, \varepsilon) \times \Sigma$ within $V$ such that each $\Sigma_t = \{t\}\times \Sigma$ has constant outward null expansion $\theta^+(t)$. The one-sided Penrose argument above applies to each $\Sigma_t$ and we conclude $\theta^+(t) \geq 0$. Following the proof of Theorem~3.1 in \cite{Galloway2018}, using the vanishing scalar curvature of the circles $\Sigma_t$, we conclude that $\theta^+(t) = 0$ for all $t \in [0, \varepsilon)$. 

Continue the foliation maximally outward, parameterized by swept area (i.e., the parameter $t$ is equal to the area of the annulus swept out between $\Sigma_0$ and $\Sigma_t$). The maximal parameter interval $[0,T)$ has $T \leq \text{Area}(V)$. We will obtain an embedded limit circle MOTS $S$ which extends the foliation via the stability and local foliation argument above, contradicting maximality. The MOTS equation gives $|A_t| \leq \norm{K}_{C^0}$, where $A_t$ and $K$ are the second fundamental forms of $\Sigma_t$ in $V$ and $V$ in $M$, respectively.  The first variation of length yields the following bound
\[
|\Sigma_t| \leq |\Sigma_0| + \norm{K}_{C^0}\text{Area}(V).
\]
Consequently, the elementary curve version of the Andersson--Metzger compactness theorem \cite[Thm. 1.3]{LarsJanJDG} gives a subsequential immersed limit $S$ as $t \to T$, which is smooth and satisfies the MOTS equation by ODE regularity. Moreover, the local graphical convergence argument in the proof of \cite[Thm. 8.1]{LarsJanJDG} shows that the limiting image $S$ is topologically
closed\footnote{Indeed, from the proof of \cite[Thm. 8.1]{LarsJanJDG}, the uniform length bound implies that,
for each $x \in V$, there are only finitely many limiting graphical
sheets in $B_V(x,s/2)$. If $x_j \in S$ converges to $x \in V$,
then, after passing to a  subsequence, the $x_j$ lie on the same limiting graphical sheet. Hence $x$ must also lie on this sheet, and we conclude that $x \in S$. We thank Jan Metzger for helpful comments regarding this point.} in $V$ and hence compact.

To show that it is embedded, note that the nesting of the regions bounded by the leaves excludes inward self-contact of this limit (cf. the proof of \cite[Thm. 7.3]{AnderssonMetzger2009}), while outward self-contact is excluded by the curve analogue of the surgery argument in \cite[Sec. 6]{AnderssonMetzger2009} along with the outward perturbation argument of \cite[Lem. 5.2]{AnderssonMetzger2009}: the surgery produces a strictly outer trapped circle that is either nonseparating or bounds an outward region of positive genus, so the corresponding infinite cyclic cover or the covering construction above, respectively, yields a contradiction to null completeness. Thus, $S$ is an embedded circle, which is two-sided since $V$ is orientable; convergence has multiplicity one because the approximating circles are embedded in its annular neighborhood. Therefore $S$ is isotopic to the preceding leaves. Again, $\lambda_1 = 0$ for $S$, so there is a local foliation by circles of constant null expansion. The same arguments above show that these circles satisfy $\theta^+ = 0$ on the outward side of $S$ but only establish $\theta^+ \geq 0$ on the inward side. Each sufficiently nearby inward circle touches a leaf of the original foliation from its outward side at a minimum of the original leaf parameter. Since the original leaf has $\theta^+ = 0$, the strong maximum principle for null expansion implies that the two circles coincide. Thus the foliation extends past the limit $S$---a contradiction.

Now assume the MOTS $\Sigma$ does not separate $V$.  Theorem 8 in \cite{GallowayLing2025} then implies that $M$ is foliated by totally geodesic null hypersurfaces. (See the remark below for some comments regarding this theorem.)  
By the proofs of Theorem \ref{Theorem: flat splitting} and Proposition \ref{Prop: splitting under flatness} in the next section, $M$ is necessarily flat (i.e. the curvature tensor vanishes)  and geodesically complete.  Theorem 3 in \cite{EhrlichGalloway} then implies that every future inextendible null geodesic contains a null cut point, and as such cannot be achronal. We arrive at the contradiction that $\gamma$ is not a null line.

We now establish the additional claim in the theorem, which follows from standard arguments. So suppose that, in addition to the assumptions from before, $(M,g)$ is timelike geodesically complete and satisfies the strong energy condition. Consider points $p_n, q_n$ with $q_n \in I^+(p_n)$ and maximizers $\gamma_n$ joining $p_n$ to $q_n$ such that $p_n$ and $q_n$ approach $-\infty$ and $+\infty$, respectively. (By this we mean, for any Cauchy surface $V'$, then there is an $N$ such that $p_n \in I^-(V')$ and $q_n \in I^+(V')$ for all $n \geq N$). Each $\gamma_n$ intersects a fixed compact Cauchy surface, so a standard limit curve argument yields, after passing to a subsequence, a causal line $\gamma$, which is necessarily timelike from the above. Then by Theorem \ref{Theorem: Lorentzian splitting}, $(M,g)$ splits isometrically as $(\R \times V, -dt^2 + h)$ and $(V,h)$ has nonnegative Gaussian curvature. By the Gauss–Bonnet theorem and the classification of compact surfaces, if $V$ is neither the two-sphere nor the projective plane, then $V$ is homeomorphic to the torus or the Klein bottle, in which case the metric $h$ is flat.
\end{proof}

\medskip

\noindent\emph{Remark.} The authors of \cite{GallowayLing2025} would like to clarify two points in the proof of \cite[Thm. 8]{GallowayLing2025}. 

First, it was not established in the proof that the family of achronal totally geodesic null hypersurfaces $\{\widetilde{N}_s\}_{s \in \R}$ in the covering spacetime $\widetilde{M}$ is pairwise disjoint. We address this now. Suppose $\widetilde{N}_s$ and $\widetilde{N}_t$ intersected at a point $\widetilde p$ for some $s < t$. Note that $\widetilde{\Sigma}_t = \widetilde{N}_t \cap \widetilde{V}$ separates $\widetilde{V}$ into an outer part $\widetilde{V}_t^+$ and an inner part $\widetilde{V}_t^-$.  Since $\widetilde{\Sigma}_s$ and $\widetilde{\Sigma}_t$ are disjoint, either $\widetilde p \in I^+(\widetilde{V})$ or $\widetilde p \in I^-(\widetilde{V})$. Assume the former; the latter is analogous. Let $\widetilde{\gamma}_s$ and $\widetilde{\gamma}_t$ be null generators on $\widetilde{N}_s$ and $\widetilde{N}_t$ from $\widetilde{\Sigma}_s$ and $\widetilde{\Sigma}_t$ to $\widetilde p$, respectively. Let $\pi \colon \widetilde{M} \to \widetilde{V}$ denote the projection along the integral curves of a timelike vector field. Consider the projected curve $\pi \circ \widetilde{\gamma}_t$. It starts off in $\widetilde{V}^+_t$. If it left this region, then it would have to meet $\widetilde{\Sigma}_t$ which contradicts achronality of $\widetilde{N}_t$. Therefore $\pi(\widetilde p) \in \widetilde{V}^+_t$. Now consider the projected curve $\pi \circ \widetilde{\gamma}_s$. It starts off in $\widetilde{V}^-_t$ and ends at $\pi(\widetilde p)$. Thus, it must intersect $\widetilde{\Sigma}_t$ yielding another achronality violation of $\widetilde{N}_t$. With pairwise disjointness established, the null normal exponential map is a diffeomorphism onto its image, giving a smooth foliation of that image. A connectedness argument shows that it covers $\widetilde{M}$. 

 Second, the foliation $\{\widetilde{N}_s\}_{s \in \R}$ of $\widetilde{M}$ descends to a foliation of  $M$, since the covering is normal and the foliation is invariant under deck transformations.

\section{A $2+1$ splitting theorem for null foliated spacetimes}\label{sec: null foliation}

In this section, we establish a splitting theorem for $3$-dimensional globally hyperbolic, null geodesically complete spacetimes with compact Cauchy surfaces which satisfy the DEC and which are foliated by totally geodesic null hypersurfaces. The latter condition implies that the spacetime belongs to the Kundt class \cite{MehidiZeghibKundt}, and we will show that the former assumptions necessarily lead to flatness in this class. This suffices by Proposition \ref{Prop: splitting under flatness} below.

\subsection{The main result}
\label{subsec: main result}

Our main result is the following:

\begin{Theorem}[$3$-dimensional cosmological splitting under null foliations]
\label{Theorem: flat splitting}
Let $(M,g)$ be a globally hyperbolic, null geodesically complete spacetime of dimension $3$ obeying the DEC. Assume that $(M,g)$ contains a smooth compact spacelike Cauchy surface and admits a smooth foliation by totally geodesic null hypersurfaces. Then
\begin{equation}
(M,g)\cong(\mathbb R\times V,-dt^2+h),
\end{equation}
where $(V,h)$ is a compact flat Riemannian surface.
In particular, $V$ is diffeomorphic to a torus or a
Klein bottle.
\end{Theorem}

\noindent \emph{Remark.} Let us note that Theorem \ref{Theorem: flat splitting} is false in higher dimensions. Indeed, Nariai spacetime $dS_2 \times S^{d-2}$ provides a counterexample in $d \geq 4$ spacetime dimensions. Observe however that Nariai spacetime does not obey the SEC, so we do not rule out that a variant of Theorem \ref{Theorem: flat splitting} could yield isometric splitting in higher dimensions assuming additionally the SEC.
\medskip\\
The subsequent subsections are dedicated to the proof of Theorem \ref{Theorem: flat splitting}. More precisely, the topological characterization of the Cauchy surfaces is obtained rather quickly (cf.\ Lemma \ref{Lemma: Topology of cauchy surface} below), and the majority of the effort is dedicated to establishing flatness of $(M,g)$. From here, the claimed isometric splitting follows from the following observation.

\begin{Proposition}[Cosmological splitting under flatness]
\label{Prop: splitting under flatness}
Let $(M,g)$ be a globally hyperbolic, null geodesically complete, flat spacetime with compact Cauchy surfaces. Then $(M,g)$ splits isometrically as $(\R \times V, -dt^2 + h)$, where $(V,h)$ is a compact flat Riemannian manifold.
\end{Proposition}
\begin{proof}
    As a consequence of \cite[Cor.\ 2]{LafuenteLopez}, flat null geodesically complete spacetimes are geodesically complete. The claimed splitting thus follows from the resolved timelike sectional curvature version of Bartnik's splitting conjecture, which was asserted in \cite[Thm.\ 5.5]{BEMG} and the necessary arguments were supplied in \cite{EhrlichGalloway}.
\end{proof}

Another application of Theorem \ref{Theorem: flat splitting} is the following corollary.

\begin{Corollary}
Let $(M,g)$ be a globally hyperbolic spacetime of dimension $3$ with compact Cauchy surfaces which is null geodesically complete and obeys the DEC. If $M$ admits a null Killing vector field, then $(M,g)$ splits as in the conclusion of Theorem \ref{Theorem: flat splitting}.
\end{Corollary}

\begin{proof}
    Let $k$ be the null Killing vector field on $M$. Since $k$ is Killing and null, we have
    \[
    g(\nabla_kk, X) = -g(\nabla_Xk,k) = -\tfrac{1}{2}Xg(k,k) = 0
    \]
    for every vector field $X$. Thus $\nabla_kk = 0$ and so its integral curves are null geodesics. Consider the smooth rank-two distribution $\mathcal{F} = k^\perp$. 
    
    We show $\mathcal{F}$ is integrable and its integral manifolds are totally geodesic. To see this, let $e$ be a local unit spacelike vector field orthogonal to $k$. Since $M$ is three-dimensional,  $\mathcal{F}$ is locally spanned by $k$ and $e$. We have
    \[
    g(\nabla_ek,e) = 0 \quad \text{ and } \quad g(\nabla_ek,k) = 0.
    \]
    These equations follow since $k$ is Killing and null, respectively. Thus $\nabla_e k \in {\rm span}\{k\}$. Together with $\nabla_k k = 0$, we see that $\nabla_X k \in  {\rm span}\{k\}$ for all $X \in \mathcal{F}$. Consequently, for any vector fields $X,Y$ tangent to $\mathcal{F}$, 
    we have
    \[
    g(\nabla_X Y, k) = - g(Y, \nabla_Xk) = 0.
    \]
    Thus $\nabla_XY \in \mathcal{F}$, and hence $[X,Y] \in \mathcal{F}$. By the Frobenius theorem, $\mathcal{F}$ determines a smooth foliation of $M$ by null hypersurfaces. Moreover, these null hypersurfaces are totally geodesic by Lemma \ref{Lemma: nabla X k in span k}. The isometric splitting then follows by Theorem \ref{Theorem: flat splitting}.
\end{proof}

\subsection{Preliminary consequences of the curvature assumptions}

The rest of this section is dedicated to the proof of Theorem \ref{Theorem: flat splitting}. We begin with some general properties of null flows and of the curvature of totally geodesic null hypersurfaces.

\begin{Lemma}[Product coordinates from the null flow]
\label{lem:null-flow-product}
Let $(M,g)$ be a globally hyperbolic spacetime, $V \subseteq M$ a compact Cauchy surface, $\mathcal F$ a smooth foliation of $M$ with immersed null hypersurfaces, $k$ be a complete future-directed affine null generator for $\mathcal F$.
Writing $\Phi_s$ for the flow of $k$, the map
\[
F:\mathbb R\times V\longrightarrow M,
\qquad F(s,x)=\Phi_s(x),
\]
is a diffeomorphism. Moreover, $F_*(\partial_s)=k$.
\end{Lemma}

\begin{proof}
Every maximal flow line of $k$ is an inextendible causal
curve, so it meets $V$ exactly once.
Together with the absence of closed causal curves,
this shows that $F$ is bijective. Since $V$ is spacelike and $k$ is null,
\[
T_xM=\mathbb R k_x\oplus T_xV.
\]
For $a\in\mathbb R$ and $v\in T_xV$,
\[
dF_{(s,x)}(a,v)
=d\Phi_s|_x(a k_x+v).
\]
The flow map $\Phi_s$ is a diffeomorphism, with inverse
$\Phi_{-s}$, so this differential is invertible.
Thus $F$ is a bijective local diffeomorphism and hence
a diffeomorphism.
Finally, the defining equation of the flow gives
$F_*(\partial_s)=k$.
\end{proof}

\begin{Lemma}[Topology of the Cauchy surface]
\label{Lemma: Topology of cauchy surface}
Let $(M,g)$ be a $3$-dimensional spacetime
with a smooth null foliation $\mathcal F$.
Let $V$ be a compact connected spacelike hypersurface.
Then $V$ is diffeomorphic to a torus or a Klein bottle.
Moreover,
\[
E:=TV\cap T\mathcal F
\]
is a smooth line subbundle of $TV$.
\end{Lemma}

\begin{proof}
Choose a smooth future-directed generator $k$ of
$(T\mathcal F)^\perp$, and let $n$ be the future-directed
unit timelike normal to $V$.
Writing $h=g|_{TV}$, the field
\[
k^\top:=k+g(k,n)n
\]
is tangent to $V$ (since it is orthogonal to $n$) and satisfies
\[
h(k^\top,k^\top)=g(k,k) + 2g(k,n) g(k,n) - g(k,n)^2 =g(k,n)^2>0.
\]
Thus $V$ admits a nowhere-zero tangent vector field.
Poincare--Hopf gives $\chi(V)=0$, so the classification
of compact connected surfaces yields the assertion. Finally, $T\mathcal F=k^\perp$ and
$g(k,X)=h(k^\top,X)$ for $X\in TV$. Hence
\[
E=(k^\top)^{\perp_h},
\]
which is a smooth line bundle.
\end{proof}

\begin{Lemma}[Invariance of the dual one-form]
\label{lem:dual-one-form-invariance}
In the setting of Lemma~\ref{lem:null-flow-product}, set $\omega:=g(k,\cdot)$ and $\alpha:=\omega|_{TV}$. Then
\[
\Phi_s^*\omega=\omega \qquad\text{for every }s\in\mathbb R.
\]
Moreover, the flow diffeomorphism $F:\mathbb R\times V\to M$
satisfies
\[
F^*\omega=\mathrm{pr}_V^*\alpha,
\]
where $\mathrm{pr}_V:\mathbb R\times V\to V$ is the projection.
\end{Lemma}

\begin{proof}
Since $k$ is affine, for every vector field $X$ we have
\[
\begin{aligned}
(\mathcal L_k\omega)(X)
&=k(g(k,X))-g(k,[k,X])\\
&=g(\nabla_k k,X)+g(k,\nabla_Xk)\\
&=\tfrac12 X(g(k,k))=0.
\end{aligned}
\]
Hence $\Phi_s^*\omega=\omega$.
Since $\omega(k)=0$, for $a\in\mathbb R$ and
$v\in T_xV$ we obtain
\[
(F^*\omega)_{(s,x)}(a,v)
=\omega_{\Phi_s(x)}(a k+d\Phi_s(v)) = (\Phi_s^*\omega)_x(v)
=\alpha_x(v),
\]
proving the second assertion.
\end{proof}

\begin{Lemma}
\label{Lemma: Ric(k,k) = 0}
Let $(M,g)$ be a spacetime and let $N$ be a smooth totally geodesic
null hypersurface in $M$. For every smooth generator $k$ of
$(TN)^\perp$, one has
\[
\operatorname{Ric}(k,k)=0.
\]
\end{Lemma}

\begin{proof}
Choose a local null frame $(k,\ell,e_1,\dots,e_{m-2})$, where $m:=\dim M$, such that the only nontrivial inner products are
\[
g(k,\ell)=-1,\qquad g(e_i,e_j)=\delta_{ij}.
\]
Note that $T N = k^\perp$, thus $e_i \in TN$. Total geodesicity gives
$\nabla_Xk\in\operatorname{span}\{k\}$ for
$X\in TN$.
Since $TN$ is involutive (as a consequence of total geodesicity), it follows that (cf.\ also Lemma \ref{Lemma: nabla X k in span k})
\[
R(X,Y)k\in\operatorname{span}\{k\}
\qquad(X,Y\in TN).
\]
In particular, $g(R(e_i,k)k,e_i)=0$ for each $i = 1,\dots,m-2$. Taking the trace in the null frame above yields
\[
\begin{aligned}
\operatorname{Ric}(k,k)
={}-g(R(k,k)k,\ell)
    -g(R(\ell,k)k,k)
   +\sum_{i=1}^{m-2}g(R(e_i,k)k,e_i)=0,
\end{aligned}
\]
using the curvature symmetries for the first two terms.
\end{proof}

\begin{Lemma}[Consequence of the dominant energy condition]
\label{Lemma: G(k,.) = -lambda g(k,.)}
Let $(M,g)$ be a spacetime obeying the DEC. Let $k\in T_pM$ be a future-directed null vector satisfying
$\operatorname{Ric}(k,k)=0$.
Then, for some $\lambda\geq0$,
\[
G(k,\cdot)=-\lambda\,g(k,\cdot).
\]
In particular,
\[
\operatorname{Ric}(k,X)=0
\qquad\text{for every }X\in k^\perp.
\]
\end{Lemma}

\begin{proof}
Define $J$ by $g(J,\cdot)=-G(k,\cdot)$.
By DEC, $J$ is future-directed causal or zero.
Moreover,
\[
g(J,k)=-G(k,k)=-\operatorname{Ric}(k,k)=0
\]
by assumption. This forces $J=\lambda k$ with $\lambda\geq0$, proving
$G(k,\cdot)=-\lambda g(k,\cdot)$.
Hence, for $X\in k^\perp$, both $G(k,X) = -\lambda g(k,X)$ and $g(k,X)$ vanish,
hence $\operatorname{Ric}(k,X)=0$.
\end{proof}

\begin{Lemma}
\label{Lemma: Ric(e,e) geq 0}
Let $(M,g)$ be a spacetime of dimension $3$ and let $p\in M$.
Let $(k,\ell,e)$ be a null frame of $T_pM$, where $k$ and $\ell$
are future-directed, $g(k,\ell)=-1$, and $e$ is a unit spacelike
vector orthogonal to both. Then
\[
\operatorname{Ric}(e,e)=2G(k,\ell).
\]
In particular, if $(M,g)$ obeys the DEC, then $\operatorname{Ric}(e,e)\geq 0$.
\end{Lemma}

\begin{proof}
Taking the metric trace in the null frame gives
\[
R_g=-2\operatorname{Ric}(k,\ell)+\operatorname{Ric}(e,e).
\]
Consequently,
\[
G(k,\ell)
=\operatorname{Ric}(k,\ell)-\frac12R_gg(k,\ell) = \Ric(k,\ell) +\frac{1}{2}R_g
=\frac12\operatorname{Ric}(e,e).
\]
DEC implies $G(k,\ell)\geq0$, proving the assertion.
\end{proof}

In the setting of Lemma~\ref{Lemma: G(k,.) = -lambda g(k,.)},
we therefore have $$\lambda=G(k,\ell)=\tfrac12\operatorname{Ric}(e,e)$$
in dimension $3$.

Let now $(M,g)$ be a spacetime of dimension $3$ and let $V$ be a smooth
spacelike hypersurface in $M$. Let $\mathcal F$ be a smooth foliation
of a neighborhood of $V$ by totally geodesic null hypersurfaces and let $k$ be a smooth future-directed generator of
$L:=(T\mathcal F)^\perp$.

By \cite[Prop.\ 4.1]{MehidiZeghibKundt} and its proof,
near each point of $V$ there are coordinates $(u,v,x)$
such that $V=\{v=0\}$, the totally geodesic null hypersurfaces are
$\{u=\mathrm{const}\}$, and
\begin{equation}
\label{eq: kundt coordinates}
g=2\,du\,dv+H(u,v,x)\,du^2
  +2W(u,v,x)\,du\,dx+q(u,x)\,dx^2,
\qquad q>0.
\end{equation}
The local generator $K:=\partial_v$ satisfies
$g(K,\cdot)=du$ and $\nabla_KK=0$. We may
write $k=aK$ for a smooth function $a>0$.

Locally choose a unit section
$e$ of $E=TV\cap T\mathcal F$, and let $\ell$
be the unique null vector satisfying
\[
g(k,\ell)=-1,\qquad g(\ell,e)=0.
\]
Define
\[
\zeta:=-g(\nabla_\ell k,e).
\]
Since $g(\nabla_\ell k,k)=0$, we have
\[
\nabla_\ell k=bk-\zeta e
\]
for some smooth local function $b$.
Total geodesicity therefore implies that, at $p\in V$,
\[
\zeta(p)=0
\iff
(\nabla_Xk)_p\in L_p
\quad\text{for every }X\in T_pM.
\]

Under $e\mapsto-e$, the scalar $\zeta$ changes sign.
Under $k\mapsto fk$ and $\ell\mapsto f^{-1}\ell$,
it remains unchanged, since
\[
\nabla_{\ell/f}(fk)
=\nabla_\ell k+\frac{\ell(f)}{f}k.
\]
Thus $\zeta^2$ defines a smooth function on $V$,
independent of the positive normalization of $k$.

By normalization invariance, we may compute $\zeta$ using
$K=\partial_v$, with $\ell$ the corresponding null vector such that $g(K,\ell)=-1$.
Then $e=q^{-1/2}\partial_x$ and $K^\flat=du$. By symmetry of the Hessian, we have $g(\nabla_XK,Y) = g(\nabla_Y K,X)$ and 
\[
(\mathcal L_Kg)(X,Y) = g(\nabla_XK,Y) + g(\nabla_Y K,X) = 2g(\nabla_X K,Y).
\]
Since $K = \partial_v$, we have $\mathcal L_K g = (\partial_v H) du^2 + (\partial_vW)(du\otimes dx + dx\otimes du)$. Using $du(\ell)=-1$, $du(e)=0$, and $dx(e) = 1/\sqrt{q}$, we obtain
\[
2g(\nabla_\ell K,e)
=(\mathcal L_Kg)(\ell,e)
=-\frac{\partial_v W}{\sqrt q}.
\]
Therefore, along $V$,
\begin{equation}
\label{eq: zeta}
\zeta=\frac{\partial_v W}{2\sqrt q},
\end{equation}
which follows from the definition of $\zeta$.

\begin{Lemma}\label{lem:zeta-ricci}
In the preceding setting, the locally defined quantities $e$ and
$\zeta$ satisfy, along $V$,
\[
\operatorname{Ric}(e,e)=2\bigl(e(\zeta)-\zeta^2\bigr).
\]
If $\operatorname{Ric}(e,e)\geq0$ along $V$, as holds in particular
under the DEC along $V$ (cf.\ Lemma \ref{Lemma: Ric(e,e) geq 0}), then
\[
e(\zeta)\geq\zeta^2.
\]
\end{Lemma}

\begin{proof}
We use the Kundt coordinates from \eqref{eq: kundt coordinates},
adapted to $V=\{v=0\}$. One can compute the component $\Ric_{xx}$ explicitly in these coordinates to get
\begin{align*}
\operatorname{Ric}_{xx}
=\partial^2_{vx}W-\frac{\partial_xq}{2q}\partial_vW-\frac12(\partial_vW)^2,
\end{align*}
where we used $\partial_v q=0$. By the preceding computation of $\zeta$, along $V$
we have (cf.\ \eqref{eq: zeta})
\[
e=q^{-1/2}\partial_x,
\qquad
\zeta=\frac{\partial_vW}{2\sqrt q}.
\]
Since $\partial_x$ is tangent to $V=\{v=0\}$,
\[
e(\zeta)
=\frac1{\sqrt q}\,
  \partial_x\left(\frac{\partial_vW}{2\sqrt q}\right)
=\frac{\partial^2_{vx}W}{2q}-\frac{\partial_xq \,\partial_vW}{4q^2}.
\]
It follows that
\[
\operatorname{Ric}(e,e)
=\frac{\operatorname{Ric}_{xx}}q
=2\bigl(e(\zeta)-\zeta^2\bigr).
\]
Thus, $\Ric(e,e) \geq 0$ implies $e(\zeta)\geq\zeta^2$.
\end{proof}

\begin{Lemma}\label{lem:zeta-vanishes}
In the setting of Lemma~\ref{lem:zeta-ricci}, suppose that $V$ is
compact and that $\operatorname{Ric}(e,e)\geq0$
along $V$. Then the locally defined quantity $\zeta$ vanishes
identically on $V$.
Consequently, for every $p\in V$,
\[
(\nabla_Xk)_p\in L_p
\qquad\text{for every }X\in T_pM.
\]
\end{Lemma}

\begin{proof}
The function $\zeta^2$ is globally defined and smooth on
$V$, so compactness gives a point $p$ where it attains
its maximum.

Suppose that $\zeta(p)^2>0$. Choose a local unit section
$e$ of $E = TV \cap T\mathcal F$ near $p$. Since $e$ is tangent to $V$,
\[
0=e(\zeta^2)(p)=2\zeta(p)e(\zeta)(p).
\]
Thus $e(\zeta)(p)=0$, contradicting
Lemma~\ref{lem:zeta-ricci}, which gives
\[
e(\zeta)(p)\geq\zeta(p)^2>0.
\]
Therefore $\zeta^2$, and hence $\zeta$, vanishes on $V$. By Lemma \ref{Lemma: nabla X k in span k}, total geodesicity gives
$\nabla_Xk\in L$ for $X\in T\mathcal F$.
The identity
\[
\nabla_\ell k=bk-\zeta e=bk
\]
supplies the remaining direction, proving the assertion.
\end{proof}

\begin{Proposition}[Parallel null line bundle]
\label{prop:parallel-null-line}
Let $(M,g)$ be a globally hyperbolic spacetime of dimension $3$
obeying the DEC. Let $V$ be a smooth compact spacelike Cauchy
surface in $M$. Let $\mathcal F$ be a smooth foliation of $M$ by
totally geodesic null hypersurfaces. Then $L:=(T\mathcal F)^\perp$
is parallel. For every smooth future-directed generator $k$ of $L$,
there is a unique smooth one-form $\beta$ on $M$ such that
\[
\nabla_Xk=\beta(X)k
\qquad\text{for every smooth vector field }X.
\]
We call $\beta$ the \emph{recurrence one-form} of $k$.
If $k$ is affine, then $\beta(k)=0$.
\end{Proposition}

\begin{proof}
Let $k$ be any smooth future-directed generator of $L$. By Bernal--Sanchez (\cite{BernalSanchez2005}) , there exists a
smooth Cauchy temporal function $\tau:M\to\mathbb R$.
Its level sets
\[
V_t:=\tau^{-1}(t)
\]
are smooth compact spacelike Cauchy hypersurfaces. Fix $p\in M$ and apply
Lemma~\ref{lem:zeta-vanishes} to
$V_{\tau(p)}$.
It follows that
\[
(\nabla_Xk)_p\in L_p
\qquad\text{for every }X\in T_pM.
\]
Since $p$ was arbitrary, $\nabla_Xk$ takes values
in $L$ everywhere on $M$.

Choose a smooth future-directed timelike vector field
$T$ on $M$. Since $g(k,T)<0$, the formula
\[
\beta(X):=\frac{g(\nabla_Xk,T)}{g(k,T)}
\]
defines a smooth one-form on $M$.
The proportionality of $\nabla_Xk$ to $k$ gives
\[
\nabla_Xk=\beta(X)k.
\]
This one-form is unique because $k$ is nowhere zero. Since every smooth local section of $L$ has the form $s=fk$, we have
\[
\nabla_Xs
=\bigl(X(f)+f\beta(X)\bigr)k\in L,
\]
which proves that $L$ is parallel. Finally, if $k$ is affine,
\[
0=\nabla_kk=\beta(k)k.
\]
Thus $\beta(k)=0$.
\end{proof}

\begin{Lemma}\label{lem:einstein-null-kernel}
Under the hypotheses of
Proposition~\ref{prop:parallel-null-line}, let $k$ be a smooth
future-directed generator of $L=(T\mathcal F)^\perp$.
Then
\[
G(k,X)=0
\qquad\text{for every vector field }X.
\]
\end{Lemma}

\begin{proof}
Fix $p\in M$ and choose a compact spacelike Cauchy
hypersurface $V'$ through $p$, as in the proof of
Proposition~\ref{prop:parallel-null-line}.
Locally along $V'$, choose $e$ and $\ell$ as in
the definition of $\zeta$. By Proposition~\ref{prop:parallel-null-line}, the bundle
$L$ is parallel. Hence $\nabla_\ell k\in L$, so
\[
\zeta=-g(\nabla_\ell k,e)=0
\]
identically near $p$ on $V'$. Since $e$ is tangent to $V'$, we have $e(\zeta)(p)=0$.
Lemma~\ref{lem:zeta-ricci} therefore gives $\operatorname{Ric}(e,e)(p)=0$.
By Lemma~\ref{Lemma: Ric(k,k) = 0}, we have
$\operatorname{Ric}(k,k)=0$ at $p$. Hence
Lemma~\ref{Lemma: G(k,.) = -lambda g(k,.)} gives
$G(k,\cdot)=-\lambda g(k,\cdot)$ at $p$.
By the observation following Lemma~\ref{Lemma: Ric(e,e) geq 0},
\[
\lambda=\tfrac12\operatorname{Ric}(e,e)=0.
\]
Thus $G(k,\cdot)=0$ at $p$.
Since $p$ was arbitrary, the conclusion holds on $M$.
\end{proof}

\begin{Lemma}\label{lem:pure-radiation}
Let $(M,g)$ be a spacetime of dimension $m \geq 3$ obeying the DEC. Let $k$ be a smooth future-directed null
vector field satisfying
\[
G(k,\cdot)=0,
\]
and set $\omega:=g(k,\cdot)$. Then there is a unique smooth function $\rho\geq0$
on $M$ such that
\[
R_g=0,
\qquad
\operatorname{Ric}=G=\rho\,\omega\otimes\omega.
\]
\end{Lemma}

\begin{proof}
Fix $p\in M$. For every future timelike vector
$X\in T_pM$, define $J_X$ by
\[
g(J_X,Y)=-G(X,Y)
\qquad\text{for every }Y\in T_pM.
\]
By the DEC, $J_X$ is future causal or zero.
Moreover, symmetry of $G$ gives
\[
g(J_X,k)=-G(X,k)=0,
\]
hence $J_X\in\operatorname{span}\{k\}$. This implies
\[
G(X,Y)=0
\qquad\text{for every }Y\in k^\perp.
\]
Since future timelike vectors span $T_pM$, this identity
holds for every $X\in T_pM$.

Choose a null frame $(k,\ell,e_1,\dots,e_{m-2})$ with $k,\ell$
future-directed, $g(k,\ell)=-1$, and $g(e_i,e_j)=\delta_{ij}$.
Since $k,e_i\in k^\perp$ and $G$ is symmetric, its only
potentially nonzero component in this frame is
\[
\rho:=G(\ell,\ell)\geq0.
\]
Using
\[
\omega(k)=\omega(e_i)=0,\qquad \omega(\ell)=-1,
\]
we obtain $G=\rho\,\omega\otimes\omega$ at $p$.

This defines a unique smooth function $\rho$ globally:
for any smooth future-directed timelike field $T$,
\[
\rho=\frac{G(T,T)}{g(k,T)^2}.
\]
Finally, taking the metric trace yields
\[
\bigg(1 - \frac{m}{2}\bigg)R_g
=\operatorname{tr}_gG
=\rho\,g(k,k)=0.
\]
Thus $R_g=0$, and the definition of $G$ gives
$\operatorname{Ric}=G$.
\end{proof}

\begin{Lemma}\label{lem:closed-recurrence}
Let $(M,g)$ be a spacetime of dimension $3$
and let $k$ be a smooth null vector field.
Set $\omega:=g(k,\cdot)$.
Suppose that a smooth one-form $\beta$ and a smooth
function $\rho$ satisfy
\[
\nabla_Xk=\beta(X)k
\qquad\text{for every vector field }X,
\]
and
\[
\operatorname{Ric}=\rho\,\omega\otimes\omega.
\]
Then
\[
R(X,Y)k=0
\qquad\text{for all vector fields }X,Y,
\]
and $\beta$ is closed.
\end{Lemma}

\begin{proof}
Taking the metric trace of the Ricci identity gives
\[
R_g=\rho\,g(k,k)=0.
\]
Let $Q$ denote the Ricci endomorphism, defined by
$g(QX,Y)=\operatorname{Ric}(X,Y)$.
Then
\[
QX=\rho\,\omega(X)k.
\]
The three-dimensional curvature identity, using
$R_g=0$ and our curvature convention, gives
\begin{equation}
\label{eq: 3d curvature}
R(X,Y)Z
=g(Y,Z)QX-g(X,Z)QY
+\operatorname{Ric}(Y,Z)X-\operatorname{Ric}(X,Z)Y.
\end{equation}
Since $\omega(k)=0$, substitution of $Z=k$ yields
\[
R(X,Y)k
=\rho\bigl(\omega(Y)\omega(X)
          -\omega(X)\omega(Y)\bigr)k
=0.
\]
On the other hand, recurrence implies
\begin{align*}
R(X,Y)k
&=\nabla_X\bigl(\beta(Y)k\bigr)
  -\nabla_Y\bigl(\beta(X)k\bigr)
  -\beta([X,Y])k\\
&=\bigl(X(\beta(Y))-Y(\beta(X))
        -\beta([X,Y])\bigr)k\\
&=d\beta(X,Y)\,k.
\end{align*}
As $k$ is nowhere zero, we conclude that $d\beta=0$.
\end{proof}

\subsection{Existence of a complete parallel null generator}

For the remainder of Section~\ref{sec: null foliation}, we suppose that
$(M,g)$ satisfies the hypotheses of Theorem~\ref{Theorem: flat splitting}.
Fix a smooth compact spacelike Cauchy surface $V$ in $M$ and a smooth
foliation $\mathcal F$ of $M$ by totally geodesic null hypersurfaces.
By Lemma~\ref{Lemma: complete affine generator}, we may fix a complete
future-directed affine generator $k$ of $L:=(T\mathcal F)^\perp$.
Let $\Phi_s$ denote its flow and set $\omega:=g(k,\cdot)$.
By Proposition~\ref{prop:parallel-null-line}, there is a unique smooth
one-form $\beta$ such that $\nabla_Xk=\beta(X)k$ for every vector field $X$.

\begin{Lemma}[Flow invariance of the recurrence one-form]
\label{lem:recurrence-flow-invariance}
The one-forms $\beta$ and $\omega$ satisfy
\[
d\beta = 0, \qquad d\omega=\beta\wedge\omega,
\qquad
\Phi_s^*\beta=\beta,
\qquad
\Phi_s^*\omega=\omega
\quad\text{for every }s\in\mathbb R.
\]
\end{Lemma}

\begin{proof}
Proposition~\ref{prop:parallel-null-line} gives
$\beta(k)=0$.
By Lemmas~\ref{lem:einstein-null-kernel}
and~\ref{lem:pure-radiation}, we have
\[
\operatorname{Ric}=\rho\,\omega\otimes\omega
\]
for a smooth function $\rho\geq0$.
Thus Lemma~\ref{lem:closed-recurrence} applies and gives
$d\beta=0$.

Metric compatibility and torsion freeness imply
\[
\begin{aligned}
d\omega(X,Y)
&=X\bigl(g(k,Y)\bigr)-Y\bigl(g(k,X)\bigr)
  -g(k,[X,Y])\\
&=g(\nabla_Xk,Y)-g(\nabla_Yk,X)\\
&=\beta(X)\omega(Y)-\beta(Y)\omega(X).
\end{aligned}
\]
Hence $d\omega=\beta\wedge\omega$. Next, for every vector field $X$,
\[
\begin{aligned}
(\mathcal L_k\beta)(X)
&=k\bigl(\beta(X)\bigr)-\beta([k,X])\\
&=d\beta(k,X)+X\bigl(\beta(k)\bigr)=0.
\end{aligned}
\]
Consequently,
\[
\frac{d}{ds}\Phi_s^*\beta
=\Phi_s^*(\mathcal L_k\beta)=0.
\]
Since $\Phi_0$ is the identity and $k$ is complete,
this proves $\Phi_s^*\beta=\beta$ for every
$s\in\mathbb R$. Finally, $\Phi_s^*\omega=\omega$ was established in
Lemma~\ref{lem:dual-one-form-invariance}.
\end{proof}

\begin{Lemma}[Variation of the metric under the null flow]
\label{lem:metric-null-flow}
The flow of $k$ satisfies
\[
\Phi_s^*g
=g+s\bigl(\beta\otimes\omega+\omega\otimes\beta\bigr)
\qquad\text{for every }s\in\mathbb R.
\]
\end{Lemma}

\begin{proof}
We first compute the Lie derivative of $g$ along $k$: Given vector fields $X,Y$ on $M$,
\[
\begin{aligned}
(\mathcal L_kg)(X,Y)
&=k\bigl(g(X,Y)\bigr)
  -g([k,X],Y)-g(X,[k,Y])\\
&=g(\nabla_Xk,Y)+g(X,\nabla_Yk)\\
&=\beta(X)\omega(Y)+\omega(X)\beta(Y).
\end{aligned}
\]
Thus
\[
\mathcal L_kg
=\beta\otimes\omega+\omega\otimes\beta.
\]
By Lemma~\ref{lem:recurrence-flow-invariance},
both $\beta$ and $\omega$ are invariant under $\Phi_s$.
Consequently,
\[
\begin{aligned}
\frac{d}{ds}\Phi_s^*g
&=\Phi_s^*(\mathcal L_kg)\\
&=(\Phi_s^*\beta)\otimes(\Phi_s^*\omega)
  +(\Phi_s^*\omega)\otimes(\Phi_s^*\beta)\\
&=\beta\otimes\omega+\omega\otimes\beta.
\end{aligned}
\]
Integrating and using $\Phi_0^*g=g$ yields the result.
\end{proof}

\begin{Lemma}[The metric in global null-flow coordinates]
\label{lem:global-null-flow-metric}
We set
\[
h:=g|_{TV},\qquad
\alpha:=\omega|_{TV},\qquad
b:=\beta|_{TV}.
\]
Then $\alpha$ is nowhere zero, and
\[
db=0,\qquad d\alpha=b\wedge\alpha.
\]
Let
\[
F:\mathbb R\times V\longrightarrow M,
\qquad F(s,x)=\Phi_s(x),
\]
be the diffeomorphism from
Lemma~\ref{lem:null-flow-product}.
Identifying tensors on $V$ with their pullbacks
under the projection $\mathbb R\times V\to V$,
we have
\[
F^*\beta=b,\qquad F^*\omega=\alpha,
\]
and
\[
F^*g
=h+ds\otimes\alpha+\alpha\otimes ds
+s\bigl(b\otimes\alpha+\alpha\otimes b\bigr).
\]
\end{Lemma}

\begin{proof}
By Lemma~\ref{lem:recurrence-flow-invariance},
the forms $\beta$ and $\omega$ are invariant under
$\Phi_s$. Moreover, $\beta(k)=\omega(k)=0$.
For $a\in\mathbb R$ and $v\in T_x V$, we therefore
obtain
\[
\begin{aligned}
(F^*\beta)_{(s,x)}(a,v)
&=\beta_{\Phi_s(x)}
  \bigl(a k_{\Phi_s(x)}+d\Phi_s|_x(v)\bigr)\\
&=(\Phi_s^*\beta)_x(v)=b_x(v).
\end{aligned}
\]
The identity $F^*\omega=\alpha$ follows in the same way,
and was already proved in
Lemma~\ref{lem:dual-one-form-invariance}. Since exterior differentiation commutes with pullback,
restricting $d\beta=0$ and $d\omega=\beta\wedge\omega$
to $V$ gives
\[
db=0,\qquad d\alpha=b\wedge\alpha.
\]
If $\alpha_x=0$, this would mean that $k_x$ is orthogonal to $T_x V$, which is impossible since $V$ is a spacelike hypersurface and $k$ is null. Thus $\alpha$ is nowhere zero. For $v,w\in T_x V$,
Lemma~\ref{lem:metric-null-flow} gives
\[
\begin{aligned}
(F^*g)_{(s,x)}\bigl((0,v),(0,w)\bigr)
= h_x(v,w)
&+s\bigl(b_x(v)\alpha_x(w)
       +\alpha_x(v)b_x(w)\bigr).
\end{aligned}
\]
Furthermore, $dF(\partial_s)=k$ and $F^*\omega=\alpha$
give
\[
(F^*g)(\partial_s,(0,v))=\alpha(v),
\qquad
(F^*g)(\partial_s,\partial_s)=0.
\]
The formula for $F^*g$ readily follows.
\end{proof}

\begin{Lemma}[Vanishing periods along positive closed curves]
\label{lem:periods-positive-curves}
Let $c:\mathbb R/\mathbb Z\to V$ be a smooth
closed curve satisfying
\[
\alpha(\dot c(t))>0
\qquad\text{for every }t.
\]
Then
\[
\int_c b=0.
\]
\end{Lemma}

\begin{proof}
Write
\[
a(t):=b(\dot c(t)),\qquad
B:=\int_0^1 a(t)\,dt=\int_c b,
\]
and suppose that $B\neq0$. For any constant $C>0$, consider the equation 
\[
\dot\tau+a(t)\tau=-C.
\]
Setting
\[
A(t):=\int_0^t a(r)\,dr,
\]
its solution with initial value $\tau(0)=\tau_0$ is
\[
\tau(t)
=e^{-A(t)}
\left(\tau_0-C\int_0^t e^{A(r)}\,dr\right).
\]
Since $A(1)=B$, the condition $\tau(1)=\tau(0)$
is equivalent to
\[
(e^B-1)\tau_0
=-C\int_0^1 e^{A(r)}\,dr.
\]
Because $B\neq0$ by assumption, this determines a unique initial
value. Since $a$ is smooth and $1$-periodic, uniqueness
for the ODE implies that the resulting solution
$\tau$ is smooth and $1$-periodic. Now choose
\[
C>
\max_{t\in\mathbb R/\mathbb Z}
\frac{h(\dot c(t),\dot c(t))}
     {2\alpha(\dot c(t))}.
\]
This maximum is finite because the denominator is
strictly positive on the compact parameter circle. The curve
\[
\Gamma(t):=F(\tau(t),c(t))
\]
is smooth and closed. By
Lemma~\ref{lem:global-null-flow-metric},
\[
\begin{aligned}
g(\dot\Gamma,\dot\Gamma)
&=(F^*g)\bigl((\dot\tau,\dot c),
              (\dot\tau,\dot c)\bigr)\\
&=h(\dot c,\dot c)
  +2\alpha(\dot c)\bigl(\dot\tau+\tau b(\dot c)\bigr)\\
&=h(\dot c,\dot c)-2C\alpha(\dot c)<0.
\end{aligned}
\]
Thus $\Gamma$ is a closed timelike curve, contradicting
global hyperbolicity. Hence $B=0$.
\end{proof}

\begin{Lemma}[Existence of a positive closed curve]
\label{lem:positive-closed-curve}
Let $S$ be a compact smooth surface without boundary,
and let $\eta$ be a smooth nowhere-zero one-form on
$S$. Then there exists a smooth closed curve
$c:\mathbb R/\mathbb Z\to S$ such that
\[
\eta(\dot c(t))>0
\qquad\text{for every }t.
\]
\end{Lemma}
\begin{proof}
Since $S$ is a surface, the smooth line distribution
$\ker \eta$ is integrable and defines a codimension-one
foliation of $S$. By the classical existence of closed transversals for
codimension-one foliations on compact manifolds
\cite[Sec.\ 4.3, p.\ 392]{Haefliger1962},
there exists a smooth closed curve
$c:\mathbb R/\mathbb Z\to S$ transverse to this foliation.
Thus $\eta(\dot c)$ is everywhere nonzero and, by
connectedness of $\mathbb R/\mathbb Z$, has constant sign.
Reversing the parametrization if necessary gives
$\eta(\dot c)>0$.
\end{proof}

\begin{Lemma}[Primitives on the universal cover]
\label{lem:forms-universal-cover}
Let
$\pi:\widetilde V\to V$ be the universal cover of $V$.
Set
\[
\widetilde\alpha:=\pi^*\alpha,
\qquad
\widetilde b:=\pi^*b.
\]
Then there exist smooth functions
$\psi,U:\widetilde V \to\mathbb R$ such that
\[
\widetilde b=d\psi,
\qquad
\widetilde\alpha=e^\psi dU.
\]
Furthermore, $dU$ is nowhere zero.
\end{Lemma}

\begin{proof}
By Lemma~\ref{lem:global-null-flow-metric},
\[
db=0,\qquad d\alpha=b\wedge\alpha.
\]
Pulling back these identities gives
\[
d\widetilde b=0,\qquad
d\widetilde\alpha
=\widetilde b\wedge\widetilde\alpha.
\]
Since $\widetilde V$ is simply connected,
there is a smooth function $\psi$ satisfying
$d\psi=\widetilde b$. Consequently,
\[
\begin{aligned}
d\bigl(e^{-\psi}\widetilde\alpha\bigr)
&=-e^{-\psi}d\psi\wedge\widetilde\alpha
  +e^{-\psi}d\widetilde\alpha\\
&=e^{-\psi}\bigl(
  -\widetilde b\wedge\widetilde\alpha
  +\widetilde b\wedge\widetilde\alpha\bigr)=0.
\end{aligned}
\]
A second application of simple connectedness gives
a smooth function $U$ such that
\[
dU=e^{-\psi}\widetilde\alpha.
\]
Since $\pi$ is a local diffeomorphism and $\alpha$
is nowhere zero, $\widetilde\alpha$ and hence $dU$
are nowhere zero.
\end{proof}

\begin{Lemma}[Deck transformations and periods]
\label{lem:deck-transformations-periods}
Let $\psi,U$ be as in
Lemma~\ref{lem:forms-universal-cover}.
For every deck transformation
$D:\widetilde V \to\widetilde V$, there exist
constants $B_D,C_D\in\mathbb R$ such that
\[
\psi\circ D=\psi+B_D,
\qquad
U\circ D=e^{-B_D}U+C_D.
\]
Moreover, if $\gamma:[0,1]\to\widetilde V$ is a
piecewise smooth path from $p$ to $D(p)$, then
$\pi\circ\gamma$ is a closed curve and
\[
B_D=\int_{\pi\circ\gamma}b.
\]
\end{Lemma}

\begin{proof}
By definition of deck transformation, $\pi\circ D=\pi$.
Since $\widetilde b=\pi^*b$ and
$\widetilde\alpha=\pi^*\alpha$, we therefore have
\[
D^*\widetilde b=\widetilde b,
\qquad
D^*\widetilde\alpha=\widetilde\alpha.
\]
Using $d\psi=\widetilde b$, we obtain
\[
d(\psi\circ D-\psi)
=D^*(d\psi)-d\psi
=D^*\widetilde b-\widetilde b
=0.
\]
Connectedness of $\widetilde V$ implies that
$\psi\circ D-\psi$ is a constant, denoted by $B_D$.

Next, using $\widetilde\alpha=e^\psi dU$, we compute
\[
e^\psi dU
=D^*\widetilde\alpha
=e^{\psi\circ D}d(U\circ D)
=e^{\psi+B_D}d(U\circ D).
\]
Thus
\[
d(U\circ D-e^{-B_D}U)=0.
\]
Again by connectedness, there is a constant $C_D$
such that
\[
U\circ D=e^{-B_D}U+C_D.
\]
Finally, if $\gamma:[0,1] \to \widetilde V$ is piecewise smooth from $p$ to $D(p)$, then $\pi(D(p))=\pi(p)$, so $\pi\circ\gamma$
is closed. Since $\widetilde b=d\psi$, we have
\[
\int_{\pi\circ\gamma}b
=\int_\gamma\pi^*b
=\int_\gamma d\psi
=\psi(D(p))-\psi(p)
=B_D.
\]
\end{proof}

\begin{Lemma}[A consequence of a nonzero period]
\label{lem:nonzero-period-function}
Suppose that $V$ is diffeomorphic to a torus.
If $B_{D_0}\neq0$ for some deck transformation $D_0$,
then there exists $f\in C^\infty(V)$ such that
\[
\alpha=df-fb.
\]
\end{Lemma}

\begin{proof}
Given a deck transformation $D$, let $a_D:=e^{-B_D}$. By Lemma~\ref{lem:deck-transformations-periods},
\[
\psi\circ D=\psi+B_D,
\qquad
U\circ D=a_DU+C_D.
\]
Since $V$ is a torus, the deck transformation group of its universal cover is isomorphic to $\pi_1(\mathbb T^2) = \mathbb Z^2$ \cite[Cor.\ 12.9]{Lee2011}.
In particular, $D_0\circ D=D\circ D_0$. The equality $U\circ D_0\circ D = U\circ D\circ D_0$ therefore yields
\begin{align}
&U \circ D_0 \circ D = U \circ D \circ D_0 \nonumber \\
\Leftrightarrow \quad &(a_{D_0}U + C_{D_0}) \circ D = (a_{D} U + C_D) \circ D_0 \nonumber \\
\Leftrightarrow \quad &a_{D_0}(a_D U + C_D) + C_{D_0} = a_D (a_{D_0}U + C_{D_0}) + C_D \nonumber\\
\Leftrightarrow \quad 
&a_{D_0}C_D+C_{D_0}
=a_DC_{D_0}+C_D. \label{eq: aD0CD}
\end{align}
Since $B_{D_0}\neq0$, we have $a_{D_0}\neq1$.
Setting
\[
u_0:=\frac{C_{D_0}}{1-a_{D_0}}.
\]
we have due to \eqref{eq: aD0CD}
\[
C_D=(1-a_D)u_0,
\qquad
U\circ D-u_0=e^{-B_D}(U-u_0)
\]
for every deck transformation $D$. Consequently, the smooth function
\[
\widetilde f:=e^\psi(U-u_0)
\]
satisfies
\[
\widetilde f\circ D
=e^{\psi+B_D}e^{-B_D}(U-u_0)
=\widetilde f.
\]
Thus, $\widetilde f$
descends to a smooth function $f:V\to\mathbb R$:
\[
\widetilde f=f\circ\pi.
\]
Using $d\psi=\widetilde b$ and
$e^\psi dU=\widetilde\alpha$, we compute
\[
d\widetilde f
=e^\psi(U-u_0)d\psi+e^\psi dU
=\widetilde f\,\widetilde b+\widetilde\alpha.
\]
Thus
\[
\pi^*(df)=\pi^*(fb+\alpha).
\]
This implies $df=fb+\alpha$, as claimed.
\end{proof}

\begin{Lemma}[Exactness of $b$ if $V \cong \mathbb T^2$]
\label{lem:exact-recurrence-torus}
If $V$ is diffeomorphic to a torus, then $b$ is exact, i.e., there
exists $\varphi_V\in C^\infty(V)$ such that
\[
b=d\varphi_V.
\]
\end{Lemma}

\begin{proof}
We use the notation of
Lemma~\ref{lem:deck-transformations-periods}.
Suppose that $B_D\neq0$ for some deck transformation $D$.
By Lemma~\ref{lem:nonzero-period-function}, there exists
$f\in C^\infty(V)$ satisfying
\[
\alpha=df-fb.
\]
By Lemma~\ref{lem:positive-closed-curve}, there is a
smooth closed curve $c:\mathbb R/\mathbb Z\to V$
such that $\alpha(\dot c)>0$.
Lemma~\ref{lem:periods-positive-curves} gives
\[
\int_c b=0.
\]
For $t\in[0,1]$, set
\[
y(t):=f(c(t)),
\qquad
a(t):=b_{c(t)}(\dot c(t)),
\qquad
A(t):=\int_0^t a(r)\,dr.
\]
Evaluating $\alpha=df-fb$ on $\dot c(t)$ gives
\[
y'(t)-a(t)y(t)=\alpha(\dot c(t)).
\]
Therefore
\begin{equation}
\label{eq: d/dt e-A(t)y(t)}
\frac{d}{dt}\bigl(e^{-A(t)}y(t)\bigr)
=e^{-A(t)}\alpha(\dot c(t))>0.
\end{equation}
On the other hand, $y(1)=y(0)$ because $c$ is closed,
while
\[
A(1)=\int_c b=0=A(0).
\]
Integrating \eqref{eq: d/dt e-A(t)y(t)} thus yields
\[
0
=\bigl[e^{-A(t)}y(t)\bigr]_0^1
=\int_0^1 e^{-A(t)}\alpha(\dot c(t))\,dt
>0,
\]
a contradiction. It follows that $B_D=0$ for every deck transformation.
By Lemma~\ref{lem:deck-transformations-periods},
$\psi$ is therefore invariant under all deck
transformations and descends to a smooth function
$\varphi_V: V\to\mathbb R$ satisfying
\[
\psi=\varphi_V\circ\pi.
\]
Since $d\psi=\pi^*b$, we obtain
\[
\pi^*(d\varphi_V)=d\psi=\pi^*b.
\]
Injectivity of pullback by the covering map $\pi$
implies $d\varphi_V=b$.
\end{proof}

\begin{Lemma}[Global exactness of the recurrence form]
\label{lem:exact-recurrence}
There exists $\varphi\in C^\infty(M)$ such that
\[
\beta=d\varphi,
\qquad
k(\varphi)=0.
\]
\end{Lemma}

\begin{proof}
We already know by Lemma \ref{Lemma: Topology of cauchy surface} that $V$ is diffeomorphic to
a torus or a Klein bottle. If $V$ is a torus, then $b$ is exact by
Lemma~\ref{lem:exact-recurrence-torus}. Suppose instead that $V$ is a Klein bottle.
Let $p:\widehat V \to V$ be its orientable
double cover, so that $\widehat V$ is a torus.
Set
\[
\widehat\alpha:=p^*\alpha,
\qquad
\widehat b:=p^*b.
\]
These forms satisfy
\[
\widehat\alpha\neq0,\qquad
d\widehat b=0,\qquad
d\widehat\alpha=\widehat b\wedge\widehat\alpha.
\]
Every smooth closed $\widehat\alpha$-positive curve
$\widehat c$ projects to an $\alpha$-positive
closed curve $p\circ\widehat c$.
Thus Lemma~\ref{lem:periods-positive-curves} gives
\[
\int_{\widehat c}\widehat b
=\int_{p\circ\widehat c}b
=0.
\]
The proof of Lemma~\ref{lem:exact-recurrence-torus}
uses only these properties of the forms on the torus.
Applying that argument to
$(\widehat V,\widehat\alpha,\widehat b)$ yields
a smooth function $\widehat\varphi$ with
$d\widehat\varphi=\widehat b$. Let now $\iota$ be the only nontrivial deck transformation
of $p$, and define
\[
\widehat\varphi_{\mathrm{av}}
:=\frac12\bigl(
\widehat\varphi+\widehat\varphi\circ\iota
\bigr).
\]
Since $p\circ\iota=p$ and $\iota^2=\mathrm{id}$,
\[
d\widehat\varphi_{\mathrm{av}}
=\frac12(\widehat b+\iota^*\widehat b)
=\widehat b,
\qquad
\widehat\varphi_{\mathrm{av}}\circ\iota
=\widehat\varphi_{\mathrm{av}}.
\]
Consequently, $\widehat\varphi_{\mathrm{av}}$
descends to a smooth function $\varphi_V$
on $V$, and $d\varphi_V=b$. In either case, choose $\varphi_V$ with
$d\varphi_V=b$.
Using the null-flow diffeomorphism
$F:\mathbb R\times V\to M$ from Lemma \ref{lem:null-flow-product}, define
\[
\varphi(F(s,x)):=\varphi_V(x).
\]
By Lemma~\ref{lem:global-null-flow-metric},
\[
F^*(d\varphi)
=\operatorname{pr}_V^*(d\varphi_V)
=\operatorname{pr}_V^*b
=F^*\beta.
\]
Since $F$ is a diffeomorphism, $d\varphi=\beta$.
Finally, $\beta(k)=0$ implies
$k(\varphi)=d\varphi(k)=0$.
\end{proof}

\begin{Lemma}[A complete parallel null generator]
\label{lem:parallel-null-generator}
Let $\varphi$ be as in
Lemma~\ref{lem:exact-recurrence} and define
\[
K:=e^{-\varphi}k.
\]
Then $K$ is a smooth, nowhere-vanishing,
future-directed null vector field satisfying
\[
\nabla K=0.
\]
Moreover, $K$ is complete.
\end{Lemma}

\begin{proof}
Since $e^{-\varphi}>0$, the vector field $K$
is nowhere vanishing, null, and future-directed.
For every smooth vector field $X$, we have
\[
\begin{aligned}
\nabla_XK
&=X(e^{-\varphi})\,k
  +e^{-\varphi}\nabla_Xk\\
&=e^{-\varphi}
  \bigl(\beta(X)-d\varphi(X)\bigr)k\\
&=0,
\end{aligned}
\]
where the last equality uses $\beta=d\varphi$
from Lemma~\ref{lem:exact-recurrence}.

Let $\Phi_s$ denote the complete flow of $k$.
Since $k(\varphi)=0$, we have
\[
\varphi(\Phi_s(p))=\varphi(p)
\]
for every $p\in M$ and $s\in\mathbb R$.
Define
\[
\Psi_t(p):=\Phi_{e^{-\varphi(p)}t}(p),
\qquad t\in\mathbb R.
\]
Then $\Psi_0(p)=p$ and
\[
\begin{aligned}
\frac{d}{dt}\Psi_t(p)
&=e^{-\varphi(p)}k_{\Psi_t(p)}\\
&=e^{-\varphi(\Psi_t(p))}k_{\Psi_t(p)}\\
&=K_{\Psi_t(p)}.
\end{aligned}
\]
Thus $\Psi$ is the flow of $K$, defined for all
real times. Hence $K$ is complete.
\end{proof}

\begin{Lemma}[An auxiliary spacelike field]
\label{lem:auxiliary-spacelike-field}
Suppose that $V$ is orientable.
Then there exist a smooth vector field $e$ and
a smooth one-form $a$ on $M$ such that
\[
g(e,e)=1,\qquad
g(e,K)=0,\qquad
\nabla_Xe=a(X)K
\]
for every smooth vector field $X$.
Moreover, $e$ is tangent to $V$ along $V$.
\end{Lemma}

\begin{proof}
Since $K$ is a positive multiple of $k$, we have
\[
E:=TV\cap K^\perp=\ker\alpha.
\]
By Lemma~\ref{lem:global-null-flow-metric},
$\alpha$ is nowhere zero.
Thus
\[
v:=\frac{\alpha^{\sharp_h}}
{|\alpha^{\sharp_h}|_h}
\]
is a smooth global unit vector field on $V$. Fix an orientation of $V$.
Let $e_0$ be the unique $h$-unit vector field
orthogonal to $v$ such that $(v,e_0)$ is positively
oriented. Then $e_0$ is a smooth global section of $E$. Let $\Psi_s$ be the complete flow of $K$.
Since $\nabla K=0$, we have
\[
(\mathcal L_Kg)(X,Y)
=g(\nabla_XK,Y)+g(X,\nabla_YK)=0.
\]
Hence every $\Psi_s$ is an isometry.
By Lemma~\ref{lem:null-flow-product}, applied to $K$,
the map
\[
F_K:\mathbb R\times V\to M,
\qquad
F_K(s,x)=\Psi_s(x),
\]
is a diffeomorphism. Define a smooth global vector field $e$ by
\[
e_{F_K(s,x)}:=d\Psi_s|_x(e_0(x)).
\]
Since $\Psi_s$ preserves both $g$ and $K$, we obtain
\[
g(e,e)=1,\qquad g(e,K)=0.
\]
Also, $e|_V=e_0$, so $e$ is tangent to $V$
along $V$. For every smooth vector field $X$,
\begin{equation}
\label{eq: nabla X e e = 0}
g(\nabla_Xe,e)=\frac12 X(g(e,e))=0,
\end{equation}
and due to the fact that $\nabla K = 0$,
\[
g(\nabla_Xe,K)
=\underbrace{X(g(e,K))}_{=0}-\underbrace{g(e,\nabla_XK)}_{=0}=0.
\]
In dimension three,
\[
K^\perp=\operatorname{span}\{K,e\}.
\]
Thus, as a consequence of \eqref{eq: nabla X e e = 0}, $\nabla_Xe$ is a multiple of $K$. These multiples define a unique smooth one-form $a$.
Indeed, choosing a smooth future-directed timelike
vector field $T$, we can write
\[
a(X)=\frac{g(\nabla_Xe,T)}{g(K,T)},
\]
whose denominator is everywhere nonzero.
Thus $\nabla_Xe=a(X)K$.
\end{proof}

\subsection{From a complete parallel null generator to flatness}

Let $K$ be the complete parallel null generator constructed in
Lemma~\ref{lem:parallel-null-generator}. We first derive a curvature
identity in the case where $V$ is orientable.

\begin{Lemma}[Curvature and the exterior derivative of $a$]
\label{lem:curvature-two-form}
Under the hypotheses and notation of
Lemma~\ref{lem:auxiliary-spacelike-field}, set
\[
\theta:=g(K,\cdot),
\qquad
\eta:=g(e,\cdot).
\]
There exists a smooth function $\widehat \rho$ on $M$ such that
\[
\operatorname{Ric}
=\widehat\rho\,\theta\otimes\theta,
\qquad
\widehat\rho\ge0,
\qquad
da=\widehat\rho\,\theta\wedge\eta.
\]
\end{Lemma}

\begin{proof}
By Lemma~\ref{lem:einstein-null-kernel}, we have $G(K,\cdot)=0$.
Applying Lemma~\ref{lem:pure-radiation} to $K$, we obtain a smooth
function $\widehat\rho\geq0$ such that $R_g=0$ and
$\operatorname{Ric}=\widehat\rho\,\theta\otimes\theta$.

By Lemma~\ref{lem:auxiliary-spacelike-field},
$\nabla_Xe=a(X)K$.
Since $\nabla K=0$, our curvature convention gives
\[
\begin{aligned}
R(X,Y)e
&=\nabla_X(a(Y)K)-\nabla_Y(a(X)K)
  -a([X,Y])K\\
&=\bigl(X(a(Y))-Y(a(X))-a([X,Y])\bigr)K\\
&=da(X,Y)K.
\end{aligned}
\]
Let $Q$ denote the Ricci endomorphism, then the Ricci formula above implies
\[
QX=\widehat\rho\,\theta(X)K.
\]
Moreover, $\theta(e)=g(K,e)=0$, so
\[
\operatorname{Ric}(X,e)=0
\]
for every $X$. Using $R_g=0$, the three-dimensional curvature
identity \eqref{eq: 3d curvature} recalled in the proof of
Lemma~\ref{lem:closed-recurrence} yields
\[
\begin{aligned}
R(X,Y)e
&=g(Y,e)QX-g(X,e)QY\\
&\quad+\operatorname{Ric}(Y,e)X
      -\operatorname{Ric}(X,e)Y\\
&=\widehat\rho\bigl(
\theta(X)\eta(Y)-\theta(Y)\eta(X)
\bigr)K\\
&=\widehat\rho\,(\theta\wedge\eta)(X,Y)K.
\end{aligned}
\]
Comparing the two expressions for $R(X,Y)e$,
and using that $K$ is nowhere zero, gives
\[
da=\widehat\rho\,\theta\wedge\eta.
\]
\end{proof}

\begin{Lemma}[Flatness]
\label{lem:flatness}
The spacetime $(M,g)$ is flat.
\end{Lemma}

\begin{proof}
We will show the claim in the case where $V$ is orientable, the non-orientable case then readily follows by going to the orientation double cover. We use the fields and forms from
Lemma~\ref{lem:curvature-two-form}.
Let $j:V\hookrightarrow M$ denote the inclusion,
and set
\[
\Omega:=j^*(\theta\wedge\eta).
\]
We first show that $\Omega$ is nowhere zero.
Since $\theta=e^{-\varphi}\omega$, we have
\[
j^*\theta=e^{-\varphi|_V}\alpha,
\]
which is nowhere zero by
Lemma~\ref{lem:global-null-flow-metric}.
By Lemma~\ref{lem:auxiliary-spacelike-field},
$e$ is tangent to $V$ along $V$, and
\[
\theta(e)=0,\qquad \eta(e)=1.
\]
Thus, for every $v\in T_p V$,
\[
\Omega_p(v,e_p)=\theta_p(v).
\]
Since $\theta_p|_{T_p V}$ is nonzero,
this proves $\Omega_p\neq0$. We endow $V$ with the orientation determined by $\Omega$.
By Lemma~\ref{lem:curvature-two-form} and
commutation of exterior differentiation with pullback,
\[
d(j^*a)
=j^*(da)
= j^*(\widehat \rho \,\theta \wedge \eta)
=(\widehat\rho|_V)\,\Omega.
\]
Because $V$ is compact and has no boundary,
Stokes' theorem gives
\[
0=\int_V d(j^*a)
=\int_V(\widehat\rho|_V)\,\Omega.
\]
Since $\widehat\rho\ge0$ and $\Omega$ is a positive
area form, continuity implies
$\widehat\rho|_V=0$.
The identity
\[
\operatorname{Ric}
=\widehat\rho\,\theta\otimes\theta
\]
therefore shows that $\operatorname{Ric}$ vanishes
at every point of $V$.

Since $K$ is parallel and complete, its flow
$\Psi_s$ consists of isometries. Hence $\operatorname{Ric}=0$ on $M$. Since $\dim M = 3$, $(M,g)$ is thus flat because of \eqref{eq: 3d curvature}.
\end{proof}

We have shown that, under the assumptions of Theorem \ref{Theorem: flat splitting}, $(M,g)$ is flat. The claims of the theorem thus follow from Lemma \ref{Lemma: Topology of cauchy surface} and Proposition \ref{Prop: splitting under flatness}.

\subsection*{Data availability statement}

There is no data associated with this manuscript.

\subsection*{Conflict of interest statement}

The authors confirm that no conflict of interest exists for this work.

\subsection*{Use of AI statement}

Discussions with ChatGPT 6 Astra were used to expedite routine computations and proofs in Section \ref{sec: null foliation}. The main ideas and final formulations are solely due to the authors, and they take full responsibility for the final contents.

\section*{Acknowledgments}

We thank Jan Metzger for helpful comments regarding the proof of Theorem \cite[Thm. 8.1]{LarsJanJDG}. The third author would like to thank the Faculty of Mathematics of the University of Vienna for the hospitality during the conference ``Curvature, Optimal Transport and Gravity" in September 2026.

This research was funded in part by the Austrian Science Fund (FWF) [Grants DOI 10.55776/EFP6, 10.55776/STA32 and 10.55776/J4913]. For open access purposes, the authors have applied a CC BY public copyright license to any author accepted manuscript version arising from this submission.
\newpage
\phantomsection
\addcontentsline{toc}{section}{References}
\bibliography{bibliography2} 
\bibliographystyle{acm}

\end{document}